\documentclass[12pt]{article}

\title{Maximal cusps are not dense}
\date{\empty}
\author{Ryo Matsuda}

\usepackage[ top=30truemm , bottom=30truemm , left=20truemm , right=20truemm]{geometry}
\usepackage{amsmath}
\usepackage{mathtools}
\usepackage{amsthm}
\usepackage{amsfonts}
\usepackage{amssymb}
\usepackage{bm}
\usepackage{color}
\usepackage[dvipdfmx]{graphicx}
\usepackage{framed}
\usepackage{empheq}
\usepackage{enumerate}
\usepackage{fancybox}
\usepackage{fancyhdr}
\usepackage{graphics}
\usepackage{mathrsfs}
\usepackage{ulem}
\usepackage{shadethm}
\usepackage{multicol}
\usepackage[all]{xy}
\usepackage{indentfirst}

\usepackage{ulem}

\newtheoremstyle{definition}
  {}
  {}
  {}
  {}
  {\bfseries}
  {.}
  { }
  {}

\theoremstyle{definition}

\newshadetheorem{de}{Definition}[section]

\newtheoremstyle{theorem}
  {}
  {12pt}
  {\itshape}
  {}
  {\bfseries}
  {.}
  { }
  {}

\theoremstyle{theorem}
\newtheorem{theo}[de]{Theorem}
\newtheorem*{copytheo}{Theorem}
\newtheorem{prop}[de]{Proposition}
\newtheorem{coro}[de]{Corollary}
\newtheorem{lem}[de]{Lemma}
\newtheorem{claim}[de]{claim}

\newtheoremstyle{answer}
  {}
  {9pt}
  {}
  {}
  {\itshape}
  {:}
  { }
  {}
  
\theoremstyle{answer}

\newtheorem*{pr}{Proof}

\newtheoremstyle{ack}
  {}
  {9pt}
  {}
  {}
  {\bfseries}
  {. }
  { }
  {}
  
\theoremstyle{ack}
\newtheorem*{ac}{Acknowledgements}
\newtheorem{Rem}[de]{Remark}


\newenvironment{proofbar}{%
   \MakeFramed {\advance\hsize-\width \FrameRestore}}%
{\endMakeFramed}

{\endMakeFramed}

\newenvironment{lembar}{
   
  \MakeFramed {\advance\hsize-\width \FrameRestore}}
{\endMakeFramed}

\definecolor{lightgray}{rgb}{0.75,0.75,0.75}
\definecolor{problemcolor}{gray}{0.3}
\definecolor{shadecolor}{gray}{0.92}

\allowdisplaybreaks[4] 

\newcommand{\Aut}{\mathrm{Aut}}
\newcommand{\AW}{\mathrm{AW}}

\newcommand{\Bel}{\mathrm{Bel}}
\newcommand{\Bers}{\mathrm{Bers}}

\newcommand{\bel}{\mathrm{bel}}

\newcommand{\Cl}{\mathrm{Cl}}

\newcommand{\dist} {\operatorname{dist}}

\newcommand{\eucli} {\mathrm{euc}}

\newcommand{\id} {\operatorname{id}}

\newcommand{\IM}{\operatorname{Im}}

\newcommand{\PSL} {PSL}

\newcommand{\short}{\mathrm{short}}
\newcommand{\leb}{\mathrm{Leb}}
\newcommand{\length}{\mathrm{length}}
\newcommand{\loc}{\mathrm{loc}}

\newcommand{\supp}{\mathrm{supp} }

\newcommand{\Teich}{\operatorname{Teich}}

\newcommand{\dr}{dr}

\newcommand{\dxdy}{\ dxdy}
\newcommand{\dxideta}{\ d\xi d\eta}

\begin{document}

\newcommand{\Addresses}{{
  \bigskip
  \footnotesize

  Ryo Matsuda, \par\nopagebreak
  \textsc{Department of Mathematics, Faculty of Science, Kyoto University, Kyoto 606-8502, Japan}\par\nopagebreak
  \textit{E-mail} : \texttt{matsuda.ryou.82c@st.kyoto-u.ac.jp}
}}


\maketitle



\begin{abstract}

	We proved that the Maximal cusp is not dense on the Bers boundary of the Teichm\"uller space of infinite-type Riemann surfaces satisfying some analytic conditions. This is a counterexample to the infinite-type case of the McMullen result for finite-type Riemann surfaces. More precisely, we showed that maximal cusps cannot approach the points on the Bers boundary corresponding to the deformation by the David map, which can be regarded as a degenerate quasiconformal map in the neighborhood of one end. In addition, to prove this, we used quasiconformal deformations in the neighborhood of a fixed end. We then proved that such a subset of the Teichm\"uller space has a manifold structure. 

\end{abstract}

\section{Introduction}

The Teichm\"uller space $\Teich(R)$ of a Riemann surface $R$ is the quotient of the equivalence relation on the set of pairs $(S, f)$ of Riemann surfaces $S$ and quasiconformal maps $f$ from $R$ to $S$. Two pairs $(S_1, f_1)$ and $(S_2, f_2)$ are Teichm\"uller equivalent if there exists a conformal map $c: S_1 \to S_2$ such that $c \circ f_1$ and $f_2$ are isotopic. If $R$ has boundaries, the dimension of $\Teich(R)$ will change depending on whether or not the isotopy fixes the boundaries.

We will briefly describe the definition of a quasiconformal map that appeared in the definition of $\Teich(R)$. A homeomorphism $f: D \to f(D)$ that preserves orientation and is defined on the plane domain $D \subset \mathbb C$ is a quasiconformal map if $f \in W_{\loc}^{1, 2} (D)$ and there exists some $k \in[0 ,1 )$ such that $| f_{\bar z}| \leq k | f_z|$ (a.e.). If we define $\mu_f : = f_{ \bar z } / f_z$, then the quasiconformla map $f$ satisfies the Beltrami equation $f_{\bar z} = \mu_f f_z$ . From the properties of quasiconformal maps, it is known that $\mu_f$ is a measurable map and satisfies $\| \mu_f\|_{L^\infty} \leq k < 1$. On the other hand, from the measurable Riemannian mapping theorem, it is known that for any $\mu \in L^\infty ( \mathbb C ) $ with $\| \mu \|_{L^\infty} < 1$ there exists a homeomorphism $f \in W_{\loc}^{1,2}(\mathbb C)$ satisfying $f_{\bar z} = \mu f_z$.
In other words, the quasiconformal maps are essentially given by $\Bel (\mathbb C): = \{ \mu \in L^\infty(\mathbb C) \mid \| \mu \|_{L^\infty} < 1 \}$. Therefore, by using the uniformization theorem and lifting to the universal covering space $\mathbb H$ of $R$, we can define the Teichm\"uller space $\Teich(\Gamma)$ of the Fuchsian group $\Gamma < \PSL(2; \mathbb R)$ such that $R \cong \mathbb H / \Gamma$. We define the space of Beltrami differentials with respect to $\Gamma$, which is described on the upper half-plane $\mathbb H$:   
	\[
		L^\infty (\mathbb H, \Gamma) 
		 : = \left\{ \mu \in L^\infty ( \mathbb H ) 
		\mid \mu \circ \gamma \frac{\bar \gamma'} {\gamma} = \mu ( \forall \gamma \in \Gamma ) \right\}, \ \ \ 
		\Bel (\mathbb H, \Gamma )  : = 
		\{ \mu \in L^\infty(\Gamma) \mid \| \mu \|_{L^\infty} < 1 \}
 	\]
If the domain is clear, write $L^\infty(\Gamma)$ and $\Bel(\Gamma)$ respectively. Also, if the Fuchsian group is a trivial group $\{ \id \}$, write $L^\infty(\mathbb H)$ and $\Bel(\mathbb H)$ respectively. 
Two elements $\mu, \nu \in \Bel ( \Gamma )$ are Teichm\"uller equivalent if there exists an element $A \in \PSL(2; \mathbb R)$ such that the self quasiconformal maps of $\mathbb H$
$f^\mu, f^\nu$ satisfies $f^\mu = A \circ f^\nu$ on the real axis $\hat { \mathbb R}$. It is known that the Teichm\"uller space $\Teich(\Gamma)$ of this Fuchsian group coincides with $\Teich(R)$.

In complex analysis, we understand a complex manifold by considering its boundary. However, it is known that the property of the boundary of a Teichm\"uller space changes greatly depending on whether $R$ is of finite or infinite type. In particular, almost nothing is known about the boundary structure when $R$ is of infinite type. For this reason, there are two types of research into the Teichm\"uller space of an infinite Riemann surface: “What conditions must be imposed to obtain the same properties as the Teichm\"uller space of a finite Riemann surface?” and “How is the Teichm\"uller space of an infinite Riemann surface different from the Teichm\"uller space of a finite Riemann surface?” This paper is classified as the latter.
Since Teichm\"uller space is the equivalence class of spaces determined by $\| \mu \| _{L^\infty} < 1$, it is thought that the Beltrami equation determined by $\| \mu \|_{L^\infty} = 1$ can be used to study its boundary. When considering the Beltrami equation on a plane domain, the conditions for the existence of solutions and the properties of solutions (for example homeomorphism, regularity...), with respect to $\| \mu \|_{L^\infty} = 1$ have been investigated by Lehto (\cite{Leh}) and David (\cite{Dav}) (\cite[Chapter 20]{AIM} systematically summarizes these). 
There are studies that use such degenerate Beltrami equations in complex dynamical systems. Based on this background, it was found that they could be applied to deformations of some Kleinian groups (\cite{main}).
In this article, we will use \cite[Theorem4.3]{main} to investigate the structure of the Bers boundary, which is the boundary of the Teichm\"uller space. The Bers boundary is obtained by embedding $\Teich(\Gamma)$ into $B(\mathbb H ^\ast,\Gamma)$ via the Bers embedding, where $B(\mathbb H ^\ast,\Gamma)$ is the Banach space of hyperbolic $L^\infty$ holomorphic quadratic differentials. 
It is also closely related to the representation space 
$\mathcal R (\Gamma) : =\{ \Gamma \to \PSL(2; \mathbb C ) \} / \text{conj}$. 
Since $\Teich(\Gamma)$ corresponds to faithful and discrete representations in $\mathcal R (\Gamma)$, we can understand that the boundary obtained when the Fuchsian group is deformed in the class of Kleinian groups is a Bers boundary. This is an important boundary that links Teichm\"uller spaces with three-dimensional hyperbolic manifold theory.
In particular, the ending lamination conjecture played a central role in this context. This conjecture was already solved by Minsky(\cite{Min}), Brock--Canary--Minsky(\cite{BCM}). 
In \cite{McM}, McMullen showed that, before the ending lamination conjecture was solved, if $R$ is of finite type, then maximal cusps are dense in the Bers boundary. 
This result is a necessary condition for the ending lamination conjecture and is a significant result. 
This is the background to this research. Next, we will state the main result:

\begin{copytheo}[Theorem \ref{main theorem}, Maximal cusps are not dense]
\begin{oframed}

	If $R$ is of infinite type with the Shiga condition and at least two ends, 
	then maximal cusps are \textbf{NOT} dense. 

\end{oframed}
\end{copytheo}

Here, an infinite Riemann surface satisfying \textit{the Shiga condition (or baounded pants decomposition)} means that there exists a pants decomposition such that the lengths $( l_j )_{j \in \mathbb N}$ of the boundaries of the decomposition satisfies
	\[
		0 < \inf_{j} l_j , \ \ \ \sup_{j} l_j < \infty. 
	\]
From the collar lemma we will discuss later, if $R$ satisfies the Shiga condition, then 
	\[
		0 < \short (R) : = \inf \{ \length_{\rho_R} ( c ) \mid c \subset R \text{: essential clpsed geodesic} \}
	\]
holds.  
The Shiga condition is a sufficient condition for the Teichm\"uller distance and the length spectrum distance to coincide and was introduced by Shiga (\cite{Shi}). There are also several studies on infinite-type Riemann surfaces and their hyperbolic geometry related to the Shiga condition; see \cite{BS} and \cite{Kin}. 
In other words, from this perspective, “an infinite Riemann surface that satisfies the Shiga condition is similar to a finite Riemann surface.” However, even if we impose such a condition (or if we impose it?), the properties of Bers boundaries differ significantly. 

In the proof of Theorem \ref{main theorem}, we used a quasiconformal deformation that only affects the neighborhood of one of the ends. We considered a subset of Teichm\"uller space that consists of such quasiconformal deformations and developed a theory of Teichm\"uller theory that focuses on the end.

\begin{copytheo}
[Theorem \ref{partial deformation at ends} \& Theorem \ref{thm: surj. par deg def}, Teichm\"uller space of partial quasiconformal deforming end]
\begin{oframed}

	The Bers projection $\mathcal B \circ \pi_T: \Bel(\mathbb H) \to B(\mathbb H^\ast)$ satisfies
		\[
			\mathcal B \circ \pi_T ( \Bel(\mathbb H, E )  ) \subset B(\mathbb H^\ast, E^\ast ). 
		\]
	And, 
	$\mathcal B \circ \pi_T : \Bel (\mathbb H, E ) \to B(\mathbb H^\ast, E^\ast ) \cap \mathcal B ( \Teich (\mathbb H ) )$
	is surjective. 

\end{oframed}
\end{copytheo}

See Section \ref{partial def around ends} for the definition of the notation. In understanding the structure of  Riemann surfaces of infinite type, 
This space $\mathcal B \circ \pi_T ( \Bel(\mathbb H, E )  ) $ can be regarded as a certain generalization of asymptotic Teichm\"uller space defined by \cite{EGL1, EGL}. 
the end (ideal boundary) is essential. In particular, quasiconformal deformations of the end are known to be crucial in infinite-dimensional Teichm\"uller theory. 
However, there are various definitions of "end". For example, when using the uniformization theorem, it is quite natural to regard
the limit set $\hat { \mathbb R }$ of the Fuchsian group as the infinite end (ideal boundary). 
In this sense, the results correspond to Theorem \ref{partial deformation at ends} and Theorem \ref{thm: surj. par deg def} can be found in \cite[Theorem 4.1]{WM}. In this paper, we give a different formulation 
so that it can be applied to the ends of Riemann surfaces as a manifold; we have clarified their structure as manifolds. For a comparison with the results in \cite[Theorem 4.1]{WM}, see Remark \ref {dif. wm and ori.} as well.


\begin{ac}
In preparing this paper, I am grateful to my supervisor, Professor Mitsuhiro Shishikura, for listening to my discussions at seminars and other occasions. 
I am also grateful to Professor Katsuhiko Matsuzaki for teaching me about \cite[Theorem 4.1]{WM} 
and for allowing me to discuss it at Waseda University. 
In addition, Professor Hiromi Ohtake has continued to hold seminars with me since my undergraduate days, 
he also helped me with various discussions about the results. 
I want to express my gratitude here. 
My friend Kento Sakai helped me remember the calculations when I finally wrote up the results,
and he helped me with various things. 
Thank you very much. 
I had the opportunity to talk about these results at many research conferences and seminars. 
All of these talks helped to improve this paper. Finally, I would like to express my deep gratitude to everyone involved.
\end{ac}

\section{Preliminaries}

\subsection{Quasihyperbolic metric}

Let $D \subset \mathbb C$ be a hyperbolic domain. We define
	\[
		\delta_D(z) : = \frac {1} { \dist_{\eucli} (z, \partial D ) }. 
	\]
The density $\delta_D(z)$  is called \textit{the quasihyperbolic metric}, and the distance determined by this is written as $d_q$ 

\begin{prop} \label{comp quasimet vs hyp met}
\begin{leftbar}
 
 	If $D = \mathbb H$, the quasihyperbolic density coincide with the hyperbolic density. 
	If $D$ is a simply connected domain, 
		\[
			\frac{1}{4} \delta_D \leq \rho_D \leq \delta_D(z), 
		\]
	where $\rho_D$ is the hyperbolic dencity of $D$. 
 
\end{leftbar}
\end{prop}

\begin{theo} [Gehring--Osgood,\ {\cite[Theorem 3]{GO}}] \label{thm: quasihyp}
\begin{oframed}
 
 	Let $D_1, D_2$ be simply connected hyperbolic domains and $f: D_1 \to D_2$ be a quasiconformal map. 
	Then, there exists a constant $c : = c (K )$ so that 
	for all $x, y \in D_1$, 
		\[
			\min \left\{ \frac{1}{c} d_{q, {D_1} } ( x, y ), \left( \frac{1}{c} d_{q, {D_1} } ( x, y ) \right)^K \right\}
			\leq 
			d_{q, {D_2} } ( f(x), f(y) ) 
			\leq
			c \max \{ d_{q, {D_1} } ( x, y ), d_{q, {D_1} } (x, y ) ^{1/K} \}, 
		\] 
 	where $d_{q, {D_j}}$ is the metric indused by the quasihyperbolic density $\rho_{D_j}$ ($j = 1, 2 $). 
\end{oframed}
\end{theo}

\subsection{Kleinian groups and Bers embedding}

A summary of Bers embedding and $B ( \Gamma )$.

\begin{de} [ ]  
	\[
		B (\mathbb H^\ast, \Gamma ) 
		: = \left\{ \varphi: \mathbb H^\ast \to \hat { \mathbb C } \middle| 
		\varphi \circ \gamma \cdot ( \gamma' ) ^2 = \varphi \ ( \forall \gamma \in \Gamma) , 
		\| \varphi \|_{B ( \Gamma )} : = \left\| \rho_{\mathbb H^\ast}^{-2} \varphi \right\|_\infty < \infty \right\}, 
	\]
where we define the hyperbolic density of the lower plane $\mathbb H^\ast$ as $\rho_{\mathbb H^\ast} : = \frac {| dz |} {- y}$. 
		
\end{de}

\begin{Rem}

	In some cases, $B(\Gamma)$ is used instead of $\mathbb H^\ast$. 
	Also, when a Kleinian group $G$ acts on a hyperbolic simply-connected domain $\Omega \subset \hat { \mathbb C }$, 
	\[
		B(\Omega, G) 
		: = \left\{ \varphi: \Omega \to \hat { \mathbb C } \mid 
		\varphi \circ \gamma \cdot ( \gamma' ) ^2 = \varphi \ ( \forall \gamma \in \Gamma) , 
		\| \rho_{\Omega}^{-2} \varphi \|_\infty < \infty \right\}. 
	\]

 \end{Rem}
	
\begin{theo} [Bers embedding (\cite{G}, 5.6, Theorem 4)] 
\begin{oframed}
		
	Let $\Gamma$ be a Fuchsian group. Then,  
	$\mathcal B : \Teich (\Gamma) \ni [\mu] \mapsto S( f_\mu ) \in B ( \Gamma )$ is a holomorphic map, and
	$\Delta_{B ( \Gamma )} ( 0 ; 1/2 ) 
	\subset \mathcal B ( \Teich ( \Gamma ) ) \subset \Cl ( \Delta_{B ( \Gamma )} ( 0 ; 3/2 ) )$.
	In particular, it is a biholomorphic map onto its image. 
	Here $f_\mu$ is the quasiconformal map with fixed $0, 1, \infty$ 
	for the Beltrami coefficient obtained by extending $\mu$ to $0$ on $\mathbb H^\ast$, and 
	\[
		S(f) : = \left( \frac {f''} {f'} \right)' ( z ) - \frac {1} {2} \left( \frac {f''} {f'} ( z ) \right)^2 .
		\]
\end{oframed}
\end{theo}

\begin{Rem}
	By using the Bers embedding, it is possible to show that $\Teich(R)$ has a complex structure. 
	In addition, it is also possible to show that this complex structure coincides 
	with the quotient complex structure of $\Bel(\Gamma)$. 
\end{Rem}

\begin{de} [Bers boundary]  

The boundary of the image of Bers embedding $\mathcal B ( \Teich ( R )$ is called \textit{Bers boundary}, and is written as  $ \partial_{\text{Bers}} \Teich(\Gamma)$. 
		
\end{de}

For each $\varphi \in B(\Gamma)$, 
then there exists a unique, locally univalent holomorphic map $W_\varphi: \mathbb H ^\ast \to \hat{\mathbb C}$ such that 
	\[
		S(W_\varphi ) = \varphi, \ \ \ W_\varphi(z) = \frac {1} {z + i} + o(z + i) \ ( \text{near} \ z = -i). 
	\] 
Then, to every an element $\gamma \in \Gamma$, there exists a unique M\"obius transformation $g$ such that 
	\[
		W_\varphi \circ \gamma = g \circ W_\varphi. 
	\]
The function $W_{\varphi}$ is called \textit{the developing map} of $\varphi$. 
We denote the mapping $\gamma \mapsto g$ by $\chi_\varphi$, called \textit{a monodromy isomorphism}. If $W_\varphi$ 
is an univalent, $\chi_\varphi(\Gamma)$ is also discrete. From Nehari's theorem (\cite[\S VI Lemma 3]{QC}),  
if $\varphi $ is contained in $\Cl(\Teich(\Gamma)) : = \Teich(\Gamma) \cup \partial_{\text{Bers}} \Teich(\Gamma)$
a function $W_\varphi$ is univalent.
Then, we classify $\varphi$ according to the properties of $\chi_\varphi(\Gamma)$. 

\begin{de} [Classification of elements in $B(\Gamma)$] \label{cal.kleiniangroup}

	Let a function $\varphi$ in $B(\Gamma)$. When the developing map of $\varphi$ is univalent, 
	the group $\chi_\varphi ( \Gamma)$ or the function $\varphi$ is called
	
		\begin{itemize}
		
			\item a \textit{topological quasiFuchsian group} 
			if the function $W_\varphi$ can be extended to a homeomorphic map on $\hat{\mathbb C}$, 
			
			\item a \textit{cusp} if there exists a hyperbolic element 
			$\gamma \in \Gamma$ such that $\chi_\varphi(\gamma)$ is a parabolic element, 
			
			\item a \textit{totally degenerate group}
			if $W_\varphi(\mathbb H^\ast) \subset \hat{\mathbb C}$ is an open and dense subset. 
			
		\end{itemize}
		
		Also, when all of the hyperbolic elements corresponding to the family of simple closed curves that give the pants
		decomposition of Riemann surface $R = \mathbb H /\Gamma$
		are mapped to the parabolic type, it is called \textit{a maximal cusp}.  
		
	\end{de} 

\begin{Rem}

 Although the word ``group'' is used, the terminology for the four Kleinian groups above will also be used for
  $\chi_\varphi$ and $\varphi$.

 \end{Rem}

Bers, Maskit investigated what Kleinian groups exist on the Bers boundary of a Teichm\"uller space in the case of 
finitely generated Fuchsian group of the first kind:

\begin{theo} [Bers\cite{B}, Maskit\cite{Mas}] 
\begin{oframed}
 
 	Let $\Gamma$ be a finitely generated Fuchsian group of the first kind. If $\chi_\varphi ( \Gamma )$ is
	a quasiFuchsian group, then $\varphi \in T ( \Gamma )$. 
	Also, if $\varphi \in \partial T ( \Gamma )$, then $\chi_\varphi ( \Gamma )$ is a cusp or 
	a totally degenerate b-group.
 
\end{oframed}
\end{theo}

\begin{Rem}

	The fact that the Fuchsian group obtained via the uniformization theorem 
	is of the first kind, which corresponds to the fact that the Riemann surface is of finite type. 

 \end{Rem}

On the other hand, using McMullen's result (Theorem \ref{McM estimate result}), 
which will be discussed later, and the degenerate Beltrami equation, the following has been proven:

\begin{theo}[{\cite[Theorem4.3]{main}}] \label {exists d f b group}
\begin{oframed}
 
 	Let $R$ be a Riemann Surface of analytically infinite type with $\short (R) > 0$ and 
	infinitely many homotopically independent essential closed geodesics 
	$\{ \gamma_n^\ast\}_{n \geq 0}$ whose lengths are bounded.  
	Then, there exists a topological quasiFuchsian group in 
	the Bers boundary of $\Teich(R)$. 
 
\end{oframed}
\end{theo}

In particular, the topological quasiFuchsian group that appears in Theorem \ref {exists d f b group} is called \textit{ a David-Fuchsian b group} because it is constructed using the David map and exists in the Bers boundary. 
 
 McMullen investigated in detail the infinitesimal behavior of Bers embedding:
 
 \begin{theo} [McMullen, {\cite[ Theorem 1.2]{McM}}] \label{McM estimate result}
\begin{oframed}
 
 	Suppose that $R$ satisfies $\short (R) > 0$. 
 	Let $[ \mu ]$ be a point in $\Teich ( R )$, 
	and $\nu$ be a unit-norm Beltrami differential in the tangent space $T_{[ \mu ]} \Teich ( R )$ at $[\mu]$. 
	If the support of $\nu$ is contained in the part of $f^\mu(R)$ of injectivity radius less than $L < 1/2$, 
	the injectivity radius is measured using the hyperbolic metric of $R$. 
	Then, the image of $\nu$ under the derivative of Bers' embedding 
	$\mathcal B: \Teich(\Gamma) \to B(\Gamma)$ has norm at most $C ( L \log 1/ L ) ^2$, 
	where the constant $C$ depends only on $\short(R)$. 	
 
\end{oframed}
\end{theo}

Using this, it was proved that there is a dense set of Maximal cusps in Bers boundaries. 

\begin{theo} [Maximal cusps are dense, {\cite[Theorem1.1]{McM}}] \label{Maximal cusps are dense}
\begin{oframed}
 
 	Let $\Gamma$ be a finitely generated Fuchsian group of the first kind, then 
	maximap cusps are dense in the Bers boundary. 
	
\end{oframed}
\end{theo}

Bers embedding has a locally inverse map. 
In particular, it is known that a very simple form of inverse map can be taken around the origin:

	\begin{theo} [Ahlfors--Weil, {\cite[Chapter VI]{QC}}] \label{AW section}
	\begin{oframed}
	
		The Bers projection $\mathcal B \circ \pi_T: \Bel ( \Gamma ) \to B(\Gamma)$ has an inverse map 
		in the neighborhood of each point of the image. 
		In particular,  in the neighborhood of the origin, it is given by %
			\[
				\AW_0: \{ \varphi \in B (\Gamma) \mid \| \varphi \|_{B(\Gamma)} < 1 /2 \} \ni \varphi 
				\mapsto -2 ( \IM z )^2 \varphi ( \bar z ) \in \Bel (\Gamma ). 
			\]	
		The map $\AW_0$ is called Ahlfors--Weil section. 
	\end{oframed}
	\end{theo}
	 
	 \begin{theo} [Ahlfors' quasiconformal reflection{\cite[Chapter IV]{QC}}] \label{Ahlfors ref}
	 \begin{oframed}
	 
	 	Let $\infty \in C \subset \hat {\mathbb C }$ be a $K$-quasicircle,
		that is, $C$ is an image of $\mathbb S^1$ by a $K$-quasiconformal map, and 
		$\Omega \cup \Omega^\ast = \hat { \mathbb C }$ be conected components.  
		There exists a sense reserving $C(K)$-quasiconformal map
		$\Lambda: \hat { \mathbb C } \to \hat { \mathbb C }$ satisfies rhe following properties: 
		real analytic, $C(K)$-Lipscitz with respect to Euclid metric, fixes $C$, and
			\[
				| \Lambda_z | < | \Lambda_{\bar z} | \leq C(K), \ \Lambda (\Omega ) = \Omega^\ast, 
				\Lambda ( \Omega^\ast ) = \Omega. 
			\]
		Such $\Lambda$ is called Ahlfors' quasiconformal reflection with respect to $C$.
	 
	 \end{oframed}
	 \end{theo}
	 
We will also summarize the relationship between Bers embedding and translation maps of Teichm\"uller space. 
In order to state the purpose in more detail, we will prepare notations. 
Let $\mu \in \Bel ( \Gamma )$.  
Using the normalized quasiconformal map $f^\mu: \mathbb H \to \mathbb H$, we define
$\Gamma^\mu : = f^\mu \circ \Gamma \circ ( f^\mu ) ^{-1}$.
Note that  $\Gamma ^\mu $ is a Fuchsian group because it acts properly discontinuous on $\mathbb H$. Then, we consider
	\[
		a_{\mu}: \Teich(\Gamma^\mu ) \ni [ \nu ] \mapsto [ \bel(f^{\nu} \circ f^{\mu} ) ] \in \Teich(\Gamma). 
	\]
Moreover, set  $\Omega = f_{\mu} ( \mathbb H )$, $\Omega^\ast = f_{\mu} ( \mathbb H^\ast )$. 

 The purpose here is to describe the relationship between the two Bers embeddings
  $\mathcal B: \Teich(\Gamma) \to B(\mathbb H^\ast, \Gamma)$, 
 $\mathcal B^{\mu}: \Teich(\Gamma^\mu) \to B(\mathbb H^\ast, \Gamma^\mu)$,  and the mapping $a_{\mu}$.
 More specifically, the purpose is to give the explicit description of the mapping $\mathcal B \circ a_{\mu} \circ ( \mathcal B^{\mu})^{-1}: B(\mathbb H, \Gamma^\mu). \to B(\mathbb H, \Gamma)$. 
	\[
		\xymatrix{
			\Teich(\Gamma^\mu) \ar[r]^-{a_{\mu}} \ar[d]^-{\mathcal B^{\mu}} & \Teich(\Gamma) \ar[d]^-{\mathcal B} \\
			B(\mathbb H^\ast, \Gamma^\mu) \ar[r]_-{ \mathcal B \circ a_{\mu} \circ ( \mathcal B^{\mu})^{-1} } 
			& B(\mathbb H^ast, \Gamma)
		}
	\]
First, we define $\Gamma_{\mu} : = f_\mu \circ \Gamma \circ ( f_{\mu} ) ^{-1}$. Note that the group is a Kleinian group.  

\begin{lem}
\begin{lembar}

	The following two maps are isometries: 
		\begin{align*}
			B(\Omega^\ast, \Gamma_\mu ) \ni \varphi 
			& \mapsto \varphi \circ f_{\mu} \cdot ( f'_{\mu} ) ^2 \in B(\mathbb H, \Gamma) \\
			B(\mathbb H^\ast, \Gamma ) \ni \varphi 
			& \mapsto ( j \varphi j ) \circ f_{\mu} \cdot ( f'_{\mu} ) ^2 \in B(\Omega, \Gamma_\mu)
		\end{align*}
\end{lembar}
\end{lem}

 \begin{theo} [{\cite[\$5.6, Theorem 2]{G}}] \label {Thm: translation map betweem Teichsp}
 \begin{oframed}
 
 	Let $\varphi \in B(\mathbb H ^\ast, \Gamma^\mu )$, $\hat { \varphi }: = ( j \varphi j ) \circ w^\mu \cdot ( w'^\mu )^2$, 
	and $\hat {\nu}: =  \rho_{\Omega}^{-2} \overline {\hat { \varphi } }$ . Then, on 
	$\Delta_{B(\mathbb H ^\ast, \Gamma^\mu )} (0 ; 1/2 )$, $\mathcal B \circ a_{\mu} \circ ( \mathcal B^{\mu})^{-1}$
	takes $\varphi$ into

		\[
			\mathcal B \circ a_{\mu} \circ ( \mathcal B^{\mu})^{-1} ( \varphi )
			= S(f_{\hat {\nu}} ) \circ f_{\mu} \circ (f'_{\mu} )^2 + S(f_\mu ), 
		\]
	where $w^\mu: \Omega \to \mathbb H$ is the conformal map that fixes $0, 1, \infty$, and $j(z) = \bar z$. 
 \end{oframed}
 \end{theo}
 
 \begin{Rem}

From the Measurable Riemann Mapping Theorem,$ w^\mu \circ f_{\mu} = f^\mu$ on   $\mathbb H$.  
From simple calculations, we can see that  $\hat { \varphi } \in B(\Omega, \Gamma_\mu )$, 
$\hat {\nu } \in \Bel ( \Omega, \Gamma_\mu )$, and 
	\begin{equation} \label {est norm for gardiner translation}
		\| \mathcal B \circ a_{\mu} \circ ( \mathcal B^{\mu})^{-1} ( \varphi ) - S(f_\mu) \|_{B(\mathbb H^\ast, \Gamma)}
		= \| S ( f_{ \hat \nu } \|_{B(\Omega, \Gamma_\mu )} = \| \varphi \|_{B(\mathbb H^\ast, \Gamma^\mu )}
		< 1 / 2. 
	\end{equation}
 \end{Rem}
 
\subsection{The collar lemma and Wolpert--Fujikawa's inequality}

For later use, we will discuss color lemmas and the Wolpert inequality.

\begin{theo}[The collor Lemma{\cite[Corollary 4.1.2]{Bus}}] \label {thecollorlem}
\begin{oframed}

	Let $R$ be a Riemann surface, and $\gamma_1, \gamma_2$ be essential simple closed geodesics. 
	If they have intersection, then
		\[
			\sinh ( \length ( \gamma_1 ) ) \cdot \sinh ( \length ( \gamma_2 ) ) > 1
		\]
\end{oframed}
\end{theo}

\begin{theo} [Wolpert, {\cite[Lemma 3.1]{W}}] \label{wolpart}
\begin{oframed}
 
 	Let $f: R \to S$ be a $K$-quasiconformal map and $\gamma$ be a closed curve  on $R$. 
	Denote by $l$ the length of the hyperbolic geodesic in the free homotopy class of $\gamma$ 
	and by $l'$ the length of the hyperbolic geodesic in the free homotopy class of $f(\gamma)$. 
	Then
			\[
				l / K \leq l' \leq lK. 
			\]

\end{oframed}
\end{theo}

Fujikawa refined Wolpert's inequality for partial quasiconformal deformations.

\begin{theo} [Fujikawa, {\cite[Lemma3.7]{Fuj}}] \label{wol. fuj. ine}
\begin{oframed}

	Let $R$ be a Riemann surface and  $c \subset R$ be an essential simple closed geodesic. 
	We set
	$d : = \dist_{\rho_{R}} ( c; E )$ for the subset $E \subset R$. 
	If $g: R \to R'$ is a K-quasiconformal map and $f|_{R \setminus E}$
	is $(1 + \varepsilon )$-quasiconformal map, then the length of  $c'$, which is homotopic to $g(c)$, satisfies
		\[
			\frac{1} {\alpha} \length(c) \leq \length ( g(c) ) \leq \alpha \length (c), 
		\] 
	where
		\[
			\alpha = K + ( 1 + \varepsilon - K ) \frac{2 \arctan ( \sinh d ) } {\pi}. 
		\]
\end{oframed}
\end{theo}

\begin{Rem}

	Suppose that we can take an essential closed geodesic $c$ that is arbitrarily far from $E$. 
	That is, for any positive number $r > 0$, 
	there exists an essential simple closed geodesic $c_r$ such that $ \dist_{\rho_{R}} ( c_r, R ) > r$. 
	In this case, $ \arctan ( \sinh d ) \to 1$ (as $d\to \infty$), so we can see that, at sufficiently large distances, 
	the way the length of the geodesic changes looks almost the same as a $(1+\varepsilon )$-quasiconformal map. 
 
 \end{Rem}

\subsection{Bergman kernel}

Here, we will review the basic calculations for the Bergman kernels of $\mathbb H, \mathbb D$ that are necessary for the later calculations.

\begin{lem}\label{rel bergman kernel d h}
\begin{lembar}

	We define three functions:
		\begin{align*}
			k_{\mathbb H}&: \mathbb H \times \mathbb H^\ast \ni ( z, \zeta ) \mapsto \frac{1}{(z - \zeta)^4} \\
			k_{\mathbb D}&: \mathbb D \times \mathbb D \ni ( z, \zeta ) \mapsto \frac{1}{(1 - z \bar\zeta)^4} \\
			G&: \hat{\mathbb C}\ni z \mapsto \frac{z - i}{z + i} \in \hat{\mathbb C}
		\end{align*}
	They have the following relation:  
		\[
			k_{\mathbb H} ( z, \zeta) = k_{\mathbb D} (G(z), G(\bar \zeta) ) ( G'(z))^2 \cdot \overline{G'(\bar \zeta)}^2. 
		\]
\end{lembar}
\end{lem}

\begin{pr}
\begin{proofbar}

	First, 
		\begin{equation*}
			G'(z) = \frac{2i} {(z + i )^2}, \ \ \ 1 - G(z)\overline{G'(\bar \zeta)} = \frac{2i(\zeta - z)}{(z+i)(\zeta - i)}. 
		\end{equation*}
	Therefore we get
		\[
			k_{\mathbb D} (G(z), G(\bar \zeta) ) ( G'(z))^2 \cdot \overline{G'(\bar \zeta)}^2 
			= \frac{1} {(2i)^4}
			\frac{(z+i)^4(\zeta-i)^4}{(\zeta - z)^4} \frac{(2i)^2}{(z+i)^2}\frac{(-2i)^4}{(\zeta-i)^4} = k_{\mathbb H} (z, \zeta). 
		\]
	\qed

\end{proofbar}
\end{pr}

\begin{lem}
\begin{lembar}

	For $A \in \Aut(\mathbb D)$, 
		\[
			k_{\mathbb D} ( Az, Aw ) (A'(z))^2 \overline{A'(w)}^2 = k_{\mathbb D}(z, w). 
		\]
\end{lembar}
\end{lem}

\begin{lem} \label{estimeate for intgral of belgman kernel over H}
\begin{lembar}

	For each $\zeta \in \mathbb H^\ast$, $| k_{\mathbb H}(z, \zeta)| $ is integrable on  $\mathbb H$. In particular, 
		\[
			\int_{\mathbb H} | k_{\mathbb H}(z, \zeta) | \dxdy = \pi \rho^2_{\mathbb H^\ast}(\zeta). 
		\]
\end{lembar}
\end{lem}

\begin{pr}
\begin{proofbar}
From the above distinctions, we get
	\begin{align*}
		\int_{\mathbb H} | k_{\mathbb H}(z, \zeta) | \dxdy 
		& \leq \int_{\mathbb H} | k_{\mathbb D} (G(z), G(\bar \zeta) )| | G'(z)|^2 \cdot |\overline{G'(\bar \zeta)}|^2 \dxdy \\
		& \quad \quad  \text{using Lemma \ref{rel bergman kernel d h}} \\
		& \int_{\mathbb D} | k_{\mathbb D} ( w, G(\bar \zeta) ) |  |\overline{G'(\bar \zeta)}|^2 \dxideta \\
		& \quad \quad ( \text{$w : = G(z), w = \xi + i \eta $ } ) \\
		& = \int_{\mathbb D} 
		| k_{\mathbb D} ( A(w), A(G(\bar \zeta)) ) | 
		|A((w)|^2 |A' (G\bar \zeta)|^2  |\overline{G'(\bar \zeta)}|^2 \dxideta \\
		&\quad \quad \text{ we take $A \in \Aut(\mathbb D)$ satisfies $A ( G (\bar \zeta) ) = 0$ } \\
		& = \int_{\mathbb D} | k_{\mathbb D} (s, 0 ) | \ dudt \cdot |A' (G\bar \zeta)|^2  |\overline{G'(\bar \zeta)}|^2 \\
		& \quad \quad \text {$s = Aw$, $s = u + it$} \\
		& = \int_{\mathbb D} \ dudt \cdot \rho^2_{\mathbb D} ( G\bar \zeta ) \cdot  |\overline{G'(\bar \zeta)}|^2 \\
		& \quad \quad \text{(The property of hyperbolic density)} \\
		& \quad \quad \text{ \quad \quad M\"obius translation $A(w) = 0$ then $\rho_{\mathbb D} (w) = | A'(w) |$} \\
		& = \pi \rho^2_{\mathbb H} ( \bar \zeta ) = \pi \rho^2_{\mathbb H^\ast} ( \zeta )
	\end{align*}
\qed
\end{proofbar}
\end{pr}

The $\pi$ appearing in the above calculation is
	\[
		\pi = \int_{\mathbb H}A \circ G  (z) \dxdy = \int_{\mathbb D} 1 \ dudt. 
	\]
If we note this, the following holds:

\begin{coro} \label{cor: core estimate}
\begin{lembar}

	Let $E \subset \mathbb H$ be a measurable set.
	For each $\zeta \in \mathbb H ^\ast$, taking $A$ and $G$  
	which apper in Lemma \ref{rel bergman kernel d h}, we get
		\[
			\int_{E} | k_{\mathbb H} ( z, \zeta ) | \dxdy = | A \circ G(E) |_{\leb} \cdot \rho^2_{\mathbb H^\ast} ( \zeta ). 
		\]
\end{lembar}
\end{coro}

\begin{lem}[Integrability for general domains] \label {lem: intable for gen case}
\begin{lembar}

	Let $C \subset \hat { \mathbb C }$ be a Jordan curve and $\Omega \cup \Omega^\ast = \hat { \mathbb C } \setminus C$
	be connected decomposition. For each $\zeta \in \Omega^\ast$, 
		\[
			\int_{\Omega} \left| \frac{1} {z - \zeta} \right|^4 \dxdy \leq 4 ^2 \pi \rho^2_{\Omega^\ast} ( \zeta ). 
		\]

\end{lembar}
\end{lem}

\begin{pr}
\begin{proofbar}

	Let $\zeta \in \Omega^\ast$. By applying the parallel translation $z \mapsto z - \zeta$ and $w \mapsto 1/w$,
	we rewrite the integral. At this time, we note that the domain is
	contained in $\hat {\mathbb C}\setminus\Delta (0;\delta_{\Omega^\ast} (\zeta))$, 	 
		\begin{align} 
			\int_{\Omega} \left| \frac{1} {z - \zeta} \right|^4 \dxdy & \leq 
			\int_{ \hat{\mathbb C} \setminus \Delta ( 0; \delta_{\Omega^\ast } ( \zeta ) )} \frac{1} {|z| ^4} \dxdy 
			\label {eq: area estimate in intable for gen case} \\
			& \leq 2 \pi \int_{\delta_{\Omega^\ast } ( \zeta ) } ^\infty \frac{1} {r^3} \dr \notag \\
			& \leq \pi ( \delta_{\Omega^\ast } ( \zeta ) ) ^{-2} . 
			\label{eq: our goal in intable for gen cass}
		\end{align}
	From Lemma \ref{comp quasimet vs hyp met}, we get
		\[
			\delta_{\Omega^\ast } ( \zeta )^{-1} \leq 4 \rho_{\Omega^\ast} ( \zeta ). 
		\]
	Therefore, 
		\[
			\int_{\Omega} \left| \frac{1} {z - \zeta} \right|^4 \dxdy \leq 4 ^2 \pi \rho^2_{\Omega^\ast} ( \zeta ). 
		\]
	 \qed
	
\end{proofbar}
\end{pr}

\begin{Rem}

Even using this method, it is possible to prove the $\leq$ in Lemma \ref{estimeate for intgral of belgman kernel over H}. 
This is because the quasihyperbolic density and the hyperbolic density coincide. Thus, we can understand that the coefficient $\pi$ in Lemma \ref{lem: intable for gen case} comes from the area in some sense. 
If we change the variable in the equation \eqref {eq: area estimate in intable for gen case} to $z \mapsto 1 / z$, we can understand that the area of a circle with radius $1 / \delta_{\Omega^\ast} ( \zeta )$ is calculated to derive. The technical part of the main result of this time is to evaluate this area. 
 \end{Rem}

\section{Maximal cusps are not dense}

\begin{theo} [Maximal cusps are not dense] \label{main theorem}
\begin{oframed}
 
 	If $R$ satisfies the Shiga condition with at least two ends, 
	then maximap cusps are not dense in $\partial_{\Bers} \Teich(R)$. 
 
\end{oframed}
\end{theo}

\begin{pr}
\begin{proofbar}

	Let $\pi : \mathbb H \to R$ be a canonical projection. Let $\varepsilon > 0$. 
	Set $m : = \short (R)$, which is greater than $0$, because
	$R$ satisfies the Shiga condition. 
 	We will prove that 
	there exists an element $\varphi \in \partial_{\Bers} \Teich(R)$ and a constant $c_0  > 0$ so that 
	for each maximal cusp $\psi \in \partial_{\Bers} \Teich(R)$, 
		\begin{equation}\label{mainestimate}
			\| \varphi - \psi \| _{B(R)} \geq c_0. 
		\end{equation}
	Since $R$ is of infinite type, there exists a sequence of essential closed geodesic $( C_n )_{n \in \mathbb N}$ such that  
		\begin{enumerate}
			\item $( C_n )$ are contained only in the subsurface $R_1$ near the one of end, 
			\item $( C_n )$ are pierwise disjoint, 
			\item $R \setminus R_1$ is connected, 
			\item $\sup_{n} (C_n) < \infty$  
		\end{enumerate}
	From Theorem \ref{exists d f b group}, there exists 
	a David--Fuchsian b group $\varphi \in \partial_{\Bers} \Teich(R)$
	that is obtained by stretching the callers of $C_n$.
	We show that this satisfies the main estimate \eqref{mainestimate}, 
	for any maximal cusp $\psi \in \partial_{\Bers} \Teich(R)$.
	
	From the construction of David--Fuchsian b group, there exists a sequence 
	$( \varphi_n )$ in $\mathcal B( \Teich(R))$ such that 
		\[
			\supp ( \bel(f_n) ) \subset R_1, \ \ \ \| \varphi_n - \varphi \|_{B(R)} \overset {n \to \infty} {\rightarrow} 0, 
		\]
	where $f_n$ is a quasiconformal extention of the developing map of $\varphi_n$. 
	Then, taking a sufficiently large $n$, set $\tilde \varphi : = \varphi_n$ satisfies
		\[
			\| \tilde \varphi - \varphi \|_{B(R)} < \varepsilon. 
		\]

	Let $\psi \in \partial_{\Bers} \Teich(R)$ be a maximal cusp and $L_j$ be the boundaries of pants which are pinched by
	$\psi$. From Theorem \ref{McM result}, there exists  $\widetilde \psi \in \mathcal B( \Teich(R))$ such that, 
		\begin{align}
			\| \widetilde \psi - \psi \|_{B(R)} &\leq \varepsilon \label{distance psi tilde psi}, \\
			\sup_{j} \length_{\tilde \psi} ( L_j ) &< m / 100 \label{lengths on tilde psi surface}.
		\end{align}
	where $ \length_{\tilde \psi}  ( L_j )$ is the length with respect to the hyperbolic metric corresponding to the $\psi$, 	
	
	Next, set $\tilde {R_1} : = \pi^{-1} ( R_1 )$, and we denote the complex conjugate of $\tilde R_1$ by  $\tilde R_1^\ast $. 
	Then we get  
		\begin{align}
			\| \varphi - \psi \|_{B(R)} &\geq \| \tilde \varphi - \tilde \psi \| - \| \tilde \varphi - \varphi \| - \| \psi - \tilde \psi \| 
			\notag \\
			& \geq \sup_{z \in \mathbb H^\ast} | ( \tilde \varphi ( z ) -\tilde \psi ( z ) )
			 \rho^{-2}_{\mathbb H^\ast} ( z ) | - 2 \varepsilon \label{first estimates} \\
			 &\geq \sup_{z \in \mathbb H^\ast \setminus (\tilde R^\ast) _r} 
			 | ( \tilde \varphi ( z ) -\tilde \psi ( z ) )
			 \rho^{-2}_{\mathbb H^\ast} ( z ) | - 2 \varepsilon \notag \\
			 &\geq  \sup_{z \in \mathbb H^\ast \setminus (\tilde R^\ast) _r}  
			 | \tilde \psi \cdot \rho^{-2}_{\mathbb H^\ast} ( z ) |
			 - 
			\sup_{z \mathbb H^\ast \setminus (\tilde R^\ast) _r}  | \tilde \varphi \cdot \rho^{-2}_{\mathbb H^\ast} ( z ) |
			- 2 \varepsilon \label{sec estimates} . 
		\end{align}
	Therefore, it is sufficient to prove that \eqref {first estimates} is grater than $1 / 4$ or 
		\begin{align}
			 \sup_{z \mathbb H^\ast \setminus (\tilde R^\ast) _r}  | \tilde \psi  \rho^{-2}_{\mathbb H^\ast} ( z ) |
			 \to c_0 > 0 \ \ \ \text{as} \ \ \ r \to \infty \label{Mcusp}.  \\ 
			 \sup_{z \mathbb H^\ast \setminus (\tilde R^\ast) _r}  | \tilde \varphi \cdot \rho^{-2}_{\mathbb H^\ast} ( z ) |
			 \to 0  \ \ \ \text{as} \ \ \ r \to \infty \label {Davidfbgroup}.
		\end{align}
The rest of the proof is given in Chapter 5.
\end{proofbar}
\end{pr}

\section{Partial deformation at end and Bers embedding} \label{partial def around ends}

\subsection{The restriction}

First, define the notation. For a measurable set  $E \subset \mathbb H$ and a positive number $r > 0$, set 
	\[
		E_r : = \{ z \in \mathbb H \mid \dist_{\rho_{\mathbb H}}(z, E ) < r \}. 
	\]
Note that E is a measurable set, so it is defined as
	\[
		\dist_{\rho_{\mathbb H}} (z, E ) = \inf \{ l > 0 \mid  | \Delta_{\rho_{\mathbb H}} ( z; l ) \cap E | _{\leb} \neq 0\}. 
	\]
We define
	\begin{align*}
		L^\infty(\mathbb H, E ) 
		&: = 
		\{ \mu \in L^\infty ( \mathbb H ) 
		\mid \limsup_{r \to \infty} \| \mu \cdot \chi_{\mathbb H \setminus E_r} \|_{L^\infty} = 0 \}, \\
		\Bel ( \mathbb H, E ) & : = \{ \mu \in L^\infty ( \mathbb H, E ) \mid \| \mu \|_{L^\infty} < 1 \}, \\
		B(\mathbb H^\ast, E^\ast ) &: = 
		\{ \varphi: \mathbb H^\ast \to \hat {\mathbb C} \mid \varphi \in B(\mathbb H),  
		\limsup_{r \to \infty} \sup_{z \in \mathbb H \setminus (E^\ast)_r} | \varphi(z) \rho_{\mathbb H^\ast} ^{-2} (z) | = 0 \}, 
		\end{align*}
where $\chi_A$ is the characteristic function of a measurable set $A$. 
	
	\begin{theo} [ ] \label{partial deformation at ends}
	\begin{oframed}
	 
	 	The Bers projection $\mathcal B \circ \pi_T: \Bel(\mathbb H) \to B(\mathbb H^\ast)$ satisfies
		\[
			\mathcal B \circ \pi_T ( \Bel(\mathbb H, E )  ) \subset B(\mathbb H^\ast, E^\ast ). 
		\]
		The same properties can also be seen in the tangent map. That is, if $\mu \in  \Bel(\mathbb H, E )$, 
		\[
			d_{[\mu]} ( \mathcal B \circ \pi_T ) : L^\infty ( \mathbb H, E ) \to B(\mathbb H^\ast, E^\ast ). 
		\]	
	\end{oframed}
	\end{theo}
	
	\begin{Rem} \label {dif. wm and ori.}
	 
	 	\cite[Theorem 4.1]{WM} also gives a result similar to Theorem\ref{partial deformation at ends}. 
		They discuss partial deformations and Bers embeddings in which the ideal boundary of the Riemann surface
		is regarded as ``end''. In this article, we discuss partial deformations and Bers embeddings i
		n which the end of the Riemann surface as a manifold is regarded as infinity. 
		In fact, if we take $\Gamma < \PSL(2; \mathbb R )$ as the first kind Fuchsian group 
		and $E$ as a set invariant under $\Gamma$, 
		then if we apply the result of \cite[Theorem 4.1]{WM} to $E$ and $\Gamma$, 
		we can only obtain a tirivial result for $\Teich(\Gamma)$. 
	 
	\end{Rem}

	\begin{lem} \label{lem: closed subsp}
	\begin{lembar}
	 
	 	$L^\infty (\mathbb H, E )$ and $B(\mathbb H^\ast, E^\ast )$ are closed subspaces in  
		$L^\infty (\mathbb H )$ and $B(\mathbb H^\ast)$ respectively. That is, they are Banach spaces.  
	 
	\end{lembar}
	\end{lem}
	
	\begin{pr}
	\begin{proofbar}
	
		Let $( \mu_n )$ be a sequence in $L^\infty (\mathbb H, E )$ such that conveges $\mu$ 
		with respect to $\|\cdot \|_{L^\infty}$. 
		Taking  $N \in \mathbb N$ satisfies
			\[
				\| \mu_n - \mu \|_{L^\infty} < \varepsilon. 
			\] 
		For $r \in (0, \infty )$, $n \geq N$, we get 
			\[
				\varepsilon > \|  \mu_n - \mu \|_{L^\infty} \geq 
				\| ( \mu_n - \mu ) \chi_{\mathbb H \setminus E_r} \|_{L^\infty}. 
			\]
		For $N$, we can take $r_N$ which satisfies
			\[
				\| \mu_N  \chi_{\mathbb H \setminus E_r} \|_{L^\infty} < \varepsilon \ \ \ ( \forall r > r_N ). 
			\]
		Therefore, we get
			\[
				\|  \mu\chi_{\mathbb H \setminus E_r} \|_{L^\infty} 
				< \varepsilon + \| \mu_N  \chi_{\mathbb H \setminus E_r} \|_{L^\infty} < 2 \varepsilon. 
			\]
		It means that $L^\infty (\mathbb H, E )$ is a closed set. The subset $L^\infty (\mathbb H, E )$ 
		is clearly a linear space, so it is a closed subspace. 
	 
	 	The second half can be proved in a similar way. 
		Let $( \varphi_n )$ be a sequence in $B(\mathbb H^\ast, E^\ast )$ such that converges $\varphi $. 
		Taking $N \in \mathbb N$ with
			\[
				\| \varphi_n - \varphi \|_{B(\mathbb H^\ast)} < \varepsilon \ \ \ ( \forall n \geq N ). 
			\]
		For  $r \in ( 0 ,\infty )$, $n \geq N$, we get 
			\[
				\varepsilon > \| \varphi_n - \varphi \|_{B(\mathbb H^\ast)} 
				\geq \sup_{z \in \mathbb H^\ast \setminus (E_r) ^\ast}
				| ( \varphi_n (z) - \varphi(z) ) \rho_{\mathbb H}^{-2} (z) |. 
			\]
		Then, for the $N$, we can take $r_N$ so that 
			\[
				 \sup_{z \in \mathbb H^\ast \setminus (E_{r}) ^\ast} | \varphi_N(z) \rho_{\mathbb H^\ast}^{-2} (z) | 
				 < \varepsilon \ \ \ ( \forall r > r_N ). 
			\]
		Therefore, we get
			\begin{align*}
				& \sup_{z \in \mathbb H^\ast \setminus (E^\ast)_r} | \varphi(z)  \rho_{\mathbb H^\ast}^{-2} (z) | 
				- 
				\sup_{z \in \mathbb H^\ast \setminus (E^\ast)_r} | \varphi_N(z) \rho_{\mathbb H^\ast}^{-2} (z) |
				 <\varepsilon \\
				& \sup_{z \in \mathbb H^\ast \setminus (E^\ast)_r} | \varphi(z)  \rho_{\mathbb H^\ast}^{-2} (z) | 
				< \sup_{z \in \mathbb H^\ast \setminus (E^\ast)_r} 
				| \varphi_N(z) \rho_{\mathbb H^\ast}^{-2} (z) | + \varepsilon. 
				< 2 \varepsilon
			\end{align*}
		That is,  $B(\mathbb H^\ast, E^\ast )$ is a closed subset.  \qed
	 
	\end{proofbar}
	\end{pr}
	
The method for proving Theorem\ref{partial deformation at ends} is the same as the method used in 
little Teichm\"uller space $\Teich_0(\Gamma)$ and $B_0(\Gamma)$,
which is the method used in \cite{EGL} to prove the correspondence by Bers embedding. In other words, 
the following will be shown in order: 
	\begin{enumerate}
		\item We show that the tangent map at the origin of the Bers embedding $d_{[0]} (\mathcal B \circ \pi_T)$ maps
		$T_0 \Bel(\mathbb H, E ) = L^\infty ( \mathbb H, E )$ to
		$T_0 B(\mathbb H, E^\ast ) = B(\mathbb H^\ast, E^\ast )$ (Proposition \ref{prop: infinitesimal claim}). 
		
		\item For $\mu \in \Bel(\mathbb H, E )$ , 
		considering the curve $g: [0, 1 ] \ni t \mapsto \mathcal B ( [t\mu] ) \in \Teich(\{ \id \} )$, 
		the derivative is $dg / dt (t) = d_{t\mu} \mathcal B( [ \mu ] )$.
		Then, we prove that$f_{t\mu} ^\ast ( \mu ) \in \Bel ( f_{t\mu} (\mathbb H^\ast) , f_{t\mu} (E^\ast ) )$
		(Corollary \ref{infinitesimal bel diff}). 
		
		\item Using the above notations, we prove 
		$d_{[t\mu]} \mathcal B( [\mu] ) \in B ( \mathbb H^\ast, E^\ast )$ (Proposition \ref{prop: gen. infini. claim}).  
	\end{enumerate}	

\begin{prop}[Infinitesimal claim of Theorem \ref{partial deformation at ends}] \label{prop: infinitesimal claim}
\begin{leftbar}
 
 	The tangent map at origin of the Bers embedding $d_{[0]} ( \mathcal B \circ \pi_T)$ satisfies 
	$d_{[0]} ( \mathcal B \circ \pi_T ) ( L^\infty ( \mathbb H, E ) ) \subset B(\mathbb H^\ast, E^\ast )$. 
 
\end{leftbar}
\end{prop}

\begin{pr}
\begin{proofbar}
 
 	Given $\mu \in \Bel(\mathbb H, E )$, $\varepsilon > 0$. Taking $r > 1$ with
		\[
			\| \mu \|_{\mathbb H \setminus E_r} < \varepsilon. 
		\]
	Note that, for $\zeta \in \mathbb H^\ast$,
		\[
			d_0 \mathcal B ( \mu ) ( \zeta ) = - \frac {6} {\pi} \int_{\mathbb H} \frac {\mu(z)} {(z - \zeta)^4} \dxdy
		\]
	( \cite[Theorem 6.11]{IT} ). 
	
	Let $\zeta \in \mathbb H^\ast \setminus E^\ast_{r}$. Separating the integral into two parts: 
		\[
			\frac {\pi} {6} | d_0 \mathcal B ( \mu ) ( \zeta ) | \leq \left( \int_{\mathbb H \setminus E_r} + \int_{E_r} \right)
			\left| \frac {\mu(z)} {(z - \zeta)^4} \right| \dxdy
		\]
	First, from Corollary\ref{cor: core estimate}, we get 
		\begin{align*}
			\int_{E_r} \left| \frac {\mu(z)} {(z - \zeta)^4} \right| \dxdy &\leq
			\int_{\mathbb H} \left| \frac {\chi_{E_r}} {(z - \zeta)^4} \right| \dxdy \\
			& = | A \circ G ( E_r ) |_{\leb} \cdot \rho^2_{\mathbb H^\ast} ( \zeta ), 
		\end{align*}
	where $G$ is defined in Lemma \ref{rel bergman kernel d h}, 
	$A \in \Aut (\mathbb D )$ satisfies $A ( G (\bar \zeta ) ) = 0$, 
	and  $z = x + iy, w = \xi + i \eta$.  
	
	Similarly, we get
		\begin{align*}
			\int_{\mathbb H \setminus E_r} \left| \frac {\mu(z)} {(z - \zeta)^4} \right| \dxdy 
			&\leq \| \mu \|_{\mathbb H \setminus E_r} \int_{\mathbb H} \left| \frac {1} {(z - \zeta)^4} \right| \dxdy \\
			& \leq \pi  \| \mu \|_{\mathbb H \setminus E_r} \cdot \rho^2_{\mathbb H^\ast} ( \zeta ). 
		\end{align*}
	Therefore
		\[
			\frac {\pi} {6} | d_0 \mathcal B ( \mu ) ( \zeta ) \rho_{\mathbb H^\ast} ( \zeta ) | \leq
			 | A \circ G ( E_r ) |_{\leb} + \pi \| \mu \|_{\mathbb H \setminus E_r}. 
		\]
	Finally, since the distance between $0$ and $A \circ G(E_r)$ is $R(r) = \frac{e^r - 1} {e^r + 1}$, 
		\begin{equation} \label{eq: area estimate}
			| A \circ G(E_r) |_{\leb} \leq \pi ( 1 - R^2(r) ). 
		\end{equation} 
	Thus, 
		\[
			\sup_{\zeta \in \mathbb H^\ast \setminus E_{2r}^\ast} 
			| d_0 \mathcal B ( \mu ) ( \zeta ) \rho_{\mathbb H^\ast} ( \zeta ) |
			\leq \frac {6} {\pi}  ( \pi ( 1 - R^2(r) ) + \pi \varepsilon. 
		\]
%
	 \qed
 
\end{proofbar}
\end{pr}

\begin{lem} \label{lem: far from far}
\begin{lembar}
 
 	Let $f: \mathbb H \to \Omega$ be a $K$-quasiconformal map. 
	Then there exists a constant $c = c(K)$ so that, for $r > c$, 
		\[
			f(E_{r/c}) \subset ( f (E) )_r. 
		\]
\end{lembar}
\end{lem}

\begin{pr}
\begin{proofbar}

	Let $c$ be a constant obtained by applying $K$ to Theorem \ref{thm: quasihyp}. 
	Take $z \in \mathbb H$ that
		\[
			r : = \dist_{\rho_{\Omega}} ( f (z), f (E) )
		\]
	is larger than $c$. Then, for $w \in E$, we get
		\[
			1 < r \leq d_{\rho, {\Omega} } ( f ( z ), f(w) ) \leq d_{q, \Omega } ( f (z), f(w) ). 
		\]
	Again, from Theorem \ref{thm: quasihyp}, it follows that  
		\begin{equation} \label{ineq: qh estimate}
			r \leq c \max \{ d_ {\rho, {\mathbb H} } (z, w ), d_{\rho, {\mathbb H} } ( z, w ) ^{1/K} \}, 
		\end{equation}
	where
		\[
			K : = \frac{1 + t\| \mu \|} {1 - t \| \mu \|}. 
		\]
	Moreover, because $r / c > 1$, 
		\[
			1 \leq \max \{ d_{\rho, {\mathbb H} } (z, w ), d_{\rho, {\mathbb H} } ( z, w ) ^{1/K} \}. 
		\]
	Indeed, if $d_{\rho, {\mathbb H} } (z, w ) \leq d_{\rho, {\mathbb H} } ( z, w ) ^{1/K}$, 
	we get $d_{\rho, {\mathbb H} } (z, w ) < 1$ since $K>1$. This is contradictory to \eqref{ineq: qh estimate}. 
	Therefore,  
		\[
			r / c \leq d_{\rho, {\mathbb H} } (z, w ). 
		\]
	That is, we get
		\[
			f (z) \not \in ( f( E) )_r \Rightarrow z \not \in E_{r/c}. 
		\]
	\qed

\end{proofbar}
\end{pr}

\begin{coro} \label{lambda exists}
\begin{lembar}
 
 	Let $\mu \in \Bel(\mathbb H, E )$. For each $t \in (0, 1 ]$,  
		\begin{equation}\label{infinitesimal bel diff}
			\lambda_t = \frac{( f_{t\mu} )_{z}} {\overline {(f_{t\mu} ) _{z}} } \cdot \frac{\mu} {1 - |t\mu|^2} \circ f_{t\mu}^{-1}
			\ ( = f_{t\mu} ^\ast ( \mu ) )
		\end{equation}
	is an element in $\Bel ( f_{t\mu} ( \mathbb H ), f_{t \mu} (E) )$ .
 
\end{lembar}
\end{coro}

\begin{pr}
\begin{proofbar}

	Note that, from the chain rule of Beltrami differential, 
		\[
			| \lambda_t | = 
			\left| \frac{( f_{t\mu} )_{z}} {\overline {(f_{t\mu} ) _{z}} } \cdot \frac{\mu} {1 - |t\mu|^2} \right| \circ f_{t\mu}^{-1}
			\leq \frac{1} {1 - t^2} | \mu | \circ f_{t\mu}^{-1}. 
		\]
	The Beltrami differntial $\mu$ satisfies
		\[
			\| \mu \chi_{\mathbb H \setminus E_r} \|_{L^\infty} \to 0 \ \ \ \text{as} \ r \to \infty. 
		\]
	Moreover, from Lemma \ref{lem: far from far}, we get
		\[
			f_{t\mu}(\mathbb H ) \setminus ( f_{t\mu} (E) )_r \subset f_{t\mu}(\mathbb H ) \setminus f_{t\mu} (E_{r/c}) 
			= f_{t\mu} ( \mathbb H \setminus E_{r/c} ). 
		\]
	Therefore, 
		\begin{align*}
			\| \mu \circ f_{t\mu}^{-1} \chi_{f_{t\mu}(\mathbb H ) \setminus ( f_{t\mu} (E) )_r} \|_{L^\infty} 
			& \leq 
			\| \mu \circ f_{t\mu}^{-1} \chi_{ f_{t\mu} ( \mathbb H \setminus E_{r/c} )} \|_{L^\infty} \\
			& = \| \mu \chi_{\mathbb H \setminus E_{r/c}} \|_{L^\infty} \to 0. 
		\end{align*}
	 \qed

\end{proofbar}
\end{pr}

\begin{lem} \label {lem: out to in}
\begin{lembar}
 
 	Let $f: \mathbb H \to \Omega$ be a $K$-quasiconformal map, 
	The constant $c$ is obtained by applying $K$ to Theorem \ref{thm: quasihyp}.
	Then $r > 4c$ implies  
		\[
			( f (E) )_{\frac{r}{4c}} \subset f(E_r). 
		\]
\end{lembar}
\end{lem}

\begin{pr}
\begin{proofbar}
 
 	Let $w \in (f(E))_{\frac{r}{4c}}$, there exists $z \in \mathbb H$ such that $w = f(z)$ and
		\[
			\dist_{\rho_\Omega} ( f(z), f(E) ) < \frac{r}{4c}. 
		\] 
	Therefore, 
		\[
			\left| \Delta_{\rho_\Omega} \left( f(z), \frac{r}{4c} \right) \cap f(E) \right|_{\leb} \neq 0. 
		\]
	We only need to show the following: 
		\begin{equation} \label{far estimate}
			\Delta_{\rho_\Omega} \left( f(z), \frac{r}{4c} \right)  \subset f ( \Delta_{\rho_\mathbb H} \left( z, r \right) ).
		\end{equation}
	 Indeed, if we prove \eqref{far estimate}, we get
		\[
			\left| \Delta_{\rho_\Omega} \left( f(z), \frac{r}{4c} \right) \cap f(E) \right|_{\leb} 
			\leq \left| f( \Delta_{\rho_\mathbb H} \left( z , r \right) \cap E ) \right|_{\leb}. 
		\]
	Moreover, since a quasiconformal map is null, 
		\[
			\left| \Delta_{\rho_\mathbb H} \left( z , r \right) \cap E  \right|_{\leb} \neq 0. 
		\]
	It means that $z \in E_r$, so  $w \in f(E_r)$. 
	
	Now, we will prove \eqref{far estimate}. 
	Let $w' = f (z' ) \in \Delta_{\rho_\Omega} \left( f(z), \frac{r}{4c} \right)$ . 
	From Theorem \ref{thm: quasihyp}, it follows that 
		\[
			 \frac{r}{4c} > d_{\rho, {\Omega} } ( f(z), w ) \geq \frac{1}{4} d_{q, {\Omega} } ( fz), f(z' ) 
			 \geq \frac{1}{4} \min \left\{ \frac{1} {c} d_{\rho, {\mathbb H} } (z, z' ), 
			 \left( \frac{1} {c} d_{\rho, {\mathbb H} } (z, z' ) \right)^K \right\}. 
		\]
	If $ d_{\rho, {\mathbb H} } (z, z' ) / c$, then $r >  d_{\rho, {\mathbb H}} (z, z' )$. It means that  
	$z' \in \Delta_{\rho_{\mathbb H^\ast}} ( z; r )$. Thus, $w' \in f( \Delta_{\rho_{\mathbb H^\ast}} ( z; r ))$. 
	Similarly, if$( d_{\rho, {\mathbb H} } (z, z' ) / c ) ^K$, then 
		\[
			\left( \frac{r}{4c} \right) ^{1/K} >  \frac{1} {c} \rho_{\mathbb H} (z, z' ). 
		\]
	Since $K \geq 1$ and $r /4c > 1$, $\left( \frac{r}{4c} \right) ^{1/K}  < r/4c$. Therefore, in both cases,
	$r >  d_{\rho, {\mathbb H} } (z, z' )$  was shown. \qed
 
\end{proofbar}
\end{pr}

Prepare the inequality that corresponds to the \eqref{eq: area estimate} that appeared 
in the proof of the proposition \ref{prop: infinitesimal claim}:

\begin{lem} \label {lem: area estimetae}
\begin{lembar}

	Let $D$ be a bounded Jordan domain such that  $C = \partial D$ is null with respect to the Lebesgue measure. 
	Let $\zeta_0 \in D$ be fixed, then the following is true:
		\[
			| D \setminus \Delta_{\rho_D} ( \zeta_0; r ) |_{\leb} \to 0 \ \ \ \text{as} \ \ \  r \to \infty. 
		\]
\end{lembar}
\end{lem}

\begin{pr}
\begin{proofbar}

	Since the hyperbolic metric is proper, $\partial  \Delta_{\rho_D} ( \zeta_0; r )$ approaches infinity. 
	In other words, it does not have an accumulation point in $D$.
	Therefore, for each$n \in \mathbb N$ , there exists $r_n$ so that 
	$r > r_0$ implies  
		\[
			\partial \Delta_{\rho_D} ( \zeta_0; r ) \subset \{ z \in \mathbb C \mid \dist (z, C ) < 1/n \}. 
		\]
	In particular, 
		\[
			D \setminus \Delta_{\rho_D} ( \zeta_0; r ) \subset \{ z \in \mathbb C \mid \dist (z, C ) < 1/n \}. 
		\]
	Moreover, we get
		\[
			\{ z \in \mathbb C \mid \dist (z, C ) < 1/n \} \to C \ \ \ \text{as} \ n \to \infty
		\]
	Since $| D |_{\leb}$ is finite, from the property of the Lebesgue measure, we get
		\[
			| \{ z \in \mathbb C \mid \dist (z, C ) < 1/n \} | _{\leb} \to | C |_{\leb} = 0 \text{as} \ n \to \infty. 
		\]
	Hence, 
		\[
			| D \setminus \Delta_{\rho_D} ( \zeta_0; r ) |_{\leb} \leq 
			| \{ z \in \mathbb C \mid \dist (z, C ) < 1/n \} | _{\leb} . 
		\]
	 \qed
			
\end{proofbar}
\end{pr}

\begin{prop} \label{prop: gen. infini. claim}
\begin{leftbar}
 
 	For each $t \in [0 , 1 ]$, 
 		\[
			d_{[t\mu]} \mathcal B ( L^\infty (\mathbb H, E ) ) \subset  B ( \mathbb H^\ast, E^\ast ). 
		\]

\end{leftbar}
\end{prop}

\begin{pr}
\begin{proofbar}

 Set $\Omega = f_{t\mu} ( \mathbb H )$, $\Omega^\ast = f_{t\mu} ( \mathbb H^\ast )$. 
Let $r > c$ and 
$\lambda_t \in \Bel ( \Omega, f_{t \mu} (E) ) $. 
For each $\zeta \in \mathbb H^\ast $ , we get  
	\[
		d_{[t\mu]} \mathcal B (\mu) = \left( - \frac{6} {\pi} 
		\int_{\Omega } \frac{\lambda_t(z)} {( z - f_{t\mu} ( \zeta ) )^4} \dxdy \right) ( f'_{t\mu} ( \zeta ) )^2
	\]
(\cite[Theorem 6.11]{IT}). By the similar calculation as Proposition \ref{prop: infinitesimal claim}, 
we divide into the following two parts: 
	\[
		\frac {\pi} {6} | d_{[t\mu]} \mathcal B (\nu_t) |
		\leq \left( \int_{ \Omega \setminus ( f_{t\mu}(E) )_r } + \int_{( f_{t\mu}(E) )_r} \right), 
		\left| \frac{\lambda_t(z)} {( z - f_{t\mu} ( \zeta ) )^4} \right| \dxdy \cdot | f'_{t\mu} ( \zeta ) | ^2
	\]
First, we get 
	\begin{align}\label{eq: estimate first 1}
		\int_{( f_{t\mu}(E) )_r} \left| \frac{\lambda_t(z)} {( z - f_{t\mu} ( \zeta ) )^4} \right|  \dxdy 
		& \leq \int_{\Omega} \left| \frac{\chi_{( f_{t\mu}(E) )_r} (z)} {( z - f_{t\mu} ( \zeta ) )^4} \right| \dxdy \\
		& \leq
		 \int_{\hat {\mathbb C } \setminus \Delta ( f_{t\mu} ( \zeta ),  \delta_{\Omega^\ast} (f_{t\mu} (\zeta ) ) )} 
		 \left| \frac{\chi_{( f_{t\mu}(E) )_r} (z)} {( z - f_{t\mu} ( \zeta ) )^4} \right| \dxdy. 
	\end{align}
Change the variable of integration, using 
	\[
		g(z) = \frac {z - f_{t\mu} ( \bar \zeta )} {z - f_{t\mu} ( \zeta )}. 
	\]
Note that 
	\[
		g^{-1} (w) = \frac{ f_{t\mu} ( \zeta ) w - f_{t\mu} ( \bar \zeta ) } {w - 1}, \ \ \ 
		( g^{-1}(w) ) ' =  \frac{ - f_{t\mu} ( \zeta ) + f_{t\mu} ( \bar \zeta ) } {( w - 1 ) ^2}. 
	\]
Set
	\[
		L : = \frac{ \delta_{\Omega^\ast} (f_{t\mu} (\zeta ) )} {|f_{t\mu} ( \zeta ) - f_{t\mu} ( \bar \zeta ) | } ( <1 ). 
	\]
Then, we get
	\[
		U: = g( \hat {\mathbb C } \setminus \Delta ( f_{t\mu} ( \zeta ),  \delta_{\Omega^\ast} (f_{t\mu} (\zeta ) ) ) ) 
		= \Delta \left( 0,  \frac{1 - L} {L} \right). 
	\]
To summarize the above,  \eqref {eq: estimate first 1}  is evaluated by
	\begin{align*}
		 & \int_{\hat {\mathbb C } \setminus \Delta ( f_{t\mu} ( \zeta ),  \delta_{\Omega^\ast} (f_{t\mu} (\zeta ) ) )} 
		 \left| \frac{\chi_{( f_{t\mu}(E) )_r} (z)} {( z - f_{t\mu} ( \zeta ) )^4} \right| \dxdy \\
		  & = 
		  \int_{U}
		  \chi_{( f_{t\mu}(E) )_r} (g^{-1}(w) ) \left| \frac{1} {f_{t\mu} ( \zeta ) - f_{t\mu} ( \bar \zeta )} \right|^2 \dxideta \\
		 &\leq 
		 \left( \frac {1} {|f_{t\mu} ( \zeta ) - f_{t\mu} ( \bar \zeta ) | } \right)^2
		  \int_{U}
		  \chi_{( f_{t\mu}(E) )_r} (g^{-1}(w) )  \dxideta \\
		  & \leq
		 \left( \frac {1} {\delta_{\Omega^\ast} ( f_{t\mu} ( \zeta )} \right)^2  \int_{g(\Omega)}
		  \chi_{( f_{t\mu}(E) )_r} (g^{-1}(w) )  \dxideta = : ( \ast )
	\end{align*}
Moreover, for each $z \in \Omega$, we get
	\[
		\delta_{\Omega^\ast} ( f_{t\mu} ( \zeta ) )
		\geq 4 | \IM \zeta | | f'_{t\mu} ( \zeta ) | 
	\]
from Koebe's 1/4 Theorem.  
Thus, 
	\begin{equation} \label{eq: estimate for first term in main prop}
		( \ast ) = \leq 16 \cdot \rho^2_{\mathbb H^\ast}( \zeta) \cdot | f'_{t\mu} ( \zeta ) |^{-2} \cdot 
		\int_{g(\Omega)} \chi_{( f_{t\mu}(E) )_r} (g^{-1}(w) )  \dxideta. 
	\end{equation}
The same applies to first-term, 
	\begin{align*}
		& \int_{ \Omega \setminus ( f_{t\mu}(E) )_r } 
		\left| \frac{\lambda_t(z)} {( z - f_{t\mu} ( \zeta ) )^4} \right| \dxdy 
		\leq \| \lambda_t \|_{L^\infty ( \Omega \setminus ( f_{t\mu}(E) )_r  )} 
		\int_{ \Omega }
		\left| \frac{1} {( z - f_{t\mu} ( \zeta ) )^4} \right| \dxdy \\
		& \leq \| \lambda_t \|_{L^\infty ( f_{t\mu}( \mathbb H) \setminus ( f_{t\mu}(E) )_r  )} 
		\cdot 16 \pi \cdot \rho^2_{\Omega^\ast}( f_{t \mu} ( \zeta ) )
		\ \ \ ( \text {Using Lemma \ref{lem: intable for gen case} } ) \\
		& \leq 16 \| \lambda_t \|_{L^\infty ( f_{t\mu}( \mathbb H) \setminus ( f_{t\mu}(E) )_r )} 
		\cdot \rho^2_{\mathbb H^\ast}( \zeta) \cdot | f'_{t\mu} ( \zeta ) |^{-2} . 
	\end{align*}
Therefore,  we get, 
	\[
		d_{[t\mu]} \mathcal B (\nu_t) ( \zeta ) \rho_{\mathbb H^\ast}^{-2} ( \zeta ) 
		\leq 16 \left ( \pi \| \lambda_t \|_{L^\infty ( f_{t\mu}( \mathbb H) \setminus ( f_{t\mu}(E) )_r )}  + 
		\int_{g(\Omega)} \chi_{( f_{t\mu}(E) )_r} (g^{-1}(w) )  \dxideta \right). 
	\]
Note that
	\[
		\| \lambda_t \|_{L^\infty ( f_{t\mu}( \mathbb H) \setminus ( f_{t\mu}(E) )_r )} \to 0 \ \ \ \text{as} \ r \to \infty
	\]
from Lemma \ref{lambda exists} and from \eqref{eq: estimate for first term in main prop}. 
Moreover, from Lemma \ref {lem: out to in}, we get
	\[
		\int_{g(\Omega)} \chi_{( f_{t\mu}(E) )_r} (g^{-1}(w) )  \dxideta
		\leq  
		\int_{g(\Omega)} \chi_{( f_{t\mu}(E_{4cr}) )} (g^{-1}(w) )  \dxideta. 
	\]
Finally, from the definition of $g$, the Jordan domain $g( \Omega )$ is bounded. Also, since a quasiconformal map is null,  
the boundary of $g(\Omega )$ is a null set. Thus, from Lemma \ref {lem: area estimetae},  
	\[
		\int_{g(\Omega)} \chi_{( f_{t\mu}(E_{4cr}) )} (g^{-1}(w) )  \dxideta \to 0 \ \ \ \text{as} \ r \to \infty. 
	\]
\qed
\end{proofbar}
\end{pr}

\begin{pr}[proof of Theorem \ref{partial deformation at ends}]
\begin{proofbar}
 
 	Follow the strategy outlined at the beginning of this section. Let  $\mu \in \Bel(\mathbb H, E )$. Considering  
		\[
			g: [0 ,1 ] \ni t \mapsto \mathcal B ([t\mu]) \in B(\mathbb H^\ast). 
		\]
	From the chain rule, the derivative of $g$ is 
		\[
			\frac{dg}{dt} = d_{[t\mu]} \mathcal B ([\mu]).
		\]
	From Proposition \ref{prop: gen. infini. claim} and 
	Lemma \ref{lem: closed subsp}, we get $g(t) \in B(\mathbb H^\ast, E^\ast)$. \qed
 
\end{proofbar}
\end{pr}

\subsection{Bers slice}

In this section, we will prove the following theorem. 
\begin{theo} \label{thm: surj. par deg def}
\begin{oframed}

	$\mathcal B \circ \pi_T : \Bel (\mathbb H) \to B(\mathbb H^\ast, E^\ast ) \cap \mathcal B ( \Teich (\mathbb H ) )$
	is holomorphic and surjective. 

\end{oframed}
\end{theo}

To prove Theorem \ref {thm: surj. par deg def}, we will show some basic things.
	
	\begin{lem} \label{lem: pull of quad. diff}
	\begin{lembar}
	
		Let $\mu \in \Bel ( \mathbb H)$.  Set $\Omega^\ast : = f_{\mu} ( \mathbb H^\ast )$. Then, the following map 
			\[
				f^\ast: B ( \Omega^\ast ) \ni \varphi \mapsto \varphi \circ f \cdot ( f' ) ^2 \in B(\mathbb H^\ast ) 
			\]
		satisfies $f^\ast ( B ( \Omega^\ast, f(E^\ast) ) ) = B( \mathbb H^\ast, E^\ast )$. 
	
	\end{lembar}
	\end{lem}
	
	\begin{pr}
	\begin{proofbar}
	
		Let $\varphi \in B ( \Omega^\ast, f(E^\ast)$. Then, we get
			\[
				f_{\mu}^\ast ( \varphi ) (w) \cdot \rho_{\Omega^\ast}^{-2} (w) = 
				\varphi ( \zeta ) 
				\cdot \rho^{-2}_{\mathbb H ^\ast} ( \zeta ), 
			\]
		where, $w = f_{\mu} ( \zeta )$. Moreover, for sufficiently large $r > 0$, applying Lemma \ref{lem: out to in}, we get 
			\[
				\sup _{\zeta \in \mathbb H^\ast \setminus E_r}\varphi ( \zeta ) 
				\cdot \rho^{-2}_{\mathbb H ^\ast} ( \zeta )
				\leq
				\sup_{w \in \Omega^\ast \setminus f_{\mu} (E^\ast )_{r / 4c}} 
				f_{\mu} ^\ast ( \varphi ) (w) \cdot \rho_{\Omega^\ast}^{-2} (w). 
			\]
		Thus $f_{\mu} ^\ast( \varphi ) \in  B( \mathbb H^\ast, E^\ast )$. 
		Considering the inverse map of $f^\ast$, we can prove the map $f_{\mu}^\ast$ is surjective.  \qed
	\end{proofbar}
	\end{pr}
	
	\begin{lem} \label{lem: pull of bel. diff}
	\begin{lembar}
	
		Let $\mu \in \Bel ( \mathbb H, E )$, $\mathbb H^\ast$. Set  $\Omega : = f ( \mathbb H)$. 
		Then, the following map 
			\[
				f_{\mu}^{\ast}: \Bel ( \Omega ) \ni \tau \mapsto \bel ( f^\tau \circ f_{\mu} ) \in \Bel ( \mathbb H )
			\]
		satisfies $f_{\mu}^\ast ( \Bel ( \Omega, f(E) ) ) = \Bel( \mathbb H, E)$.
	
	\end{lembar}
	\end{lem}
	
	\begin{pr}
	\begin{proofbar}
	
		From the chain rule of Beltrami differentials, we get 
			\[
				\bel ( f^\tau \circ f_{\mu} ) 
				= \frac{\tau \circ f_{\mu} - \theta_{\mu} \mu } {\theta_{\mu} + \overline{\mu} \cdot \tau \circ f_{\mu}}, 
			\]
		where $\theta_{\mu} : =  ( f _{\mu} )_{z} / \overline { (f_{\mu}) _{\bar z}  }$. Thus, we get 
			\[
				| \bel ( f^\tau \circ f )| \leq 
				 \frac{| \tau \circ f | + | \bel(f) | } {1 - | \overline{\bel(f)} \cdot \tau \circ f| }. 
			\]
		Taking sufficiently large $r$, $| \tau \circ f | < \varepsilon$, $| \bel(f) |. < \varepsilon$, so 
		we get $| \bel ( f^\tau \circ f )| \to 0$. 
		Considering the inverse map of $f^\ast$, we can prove the map $f_{\mu}^\ast$ is surjective.  \qed
	
	\end{proofbar}
	\end{pr}

\begin{pr}[Theorem \ref {thm: surj. par deg def}]
\begin{proofbar}

	First, consider the neighborhood of the origin. t means that, we will prove
		\[
			\{ \varphi \in B ( \mathbb H^\ast ) \mid \| \varphi \|_{B(\mathbb H^\ast )} < 1 / 2 \} 
			\cap B(\mathbb H^\ast, E^\ast )
			 \subset 
			\mathcal B \circ \pi_T (  \Bel (\mathbb H, E ) ). 
		\]
	This can be done using the Ahlfors-Weil section. Indeed, 
	For each $\varphi \in \Delta_{B(\mathbb H^\ast )} ( 0; 1/2) \cap B(\mathbb H^\ast, E^\ast )$, we define
		\[
			\AW_0 (\varphi) ( z ) = -2 ( \IM z ) ^2 \cdot \varphi ( \bar z ) \ \ \ z \in \mathbb H. 
		\]
	From Theorem \ref{AW section}, this satisfies $\mathcal B \circ \pi_T ( \AW_0 ( \varphi ) ) = \varphi$. 
	Moreover, for each $r > 0$, we get 
		\begin{align*}
			& \sup_{\zeta \in \mathbb H^\ast \setminus ( E^\ast )_r} 
			| \varphi ( \zeta ) \rho^{-2}_{\mathbb H^\ast} (\zeta) | 
			= 
			\sup_{\zeta \in \mathbb H^\ast \setminus ( E^\ast )_r} 
			| \AW_0 (\varphi) ( \bar \zeta ) | \\
			& =
			\sup_{z \in \mathbb H \setminus E_r} | \AW_0 (\varphi) ( z ) | 
			=
			\| | \AW_0 (\varphi) ( z) \chi_{\mathbb H \setminus E_r} \| _{L^\infty}, 
		\end{align*}
	so
		\[
			\| | \AW_0 (\varphi) \cdot \chi_{\mathbb H \setminus E_r} \| _{L^\infty} \to 0 \ \ \ \text{as} \ 
			r \to \infty. 
		\]
	Thus, $\AW_0 (\varphi) \in \Bel ( \mathbb H, E )$. 
	That is, $\varphi \in \mathcal B \circ \pi_T ( \Bel(\mathbb H, E ) )$. 
	
	Next, consider in a general situation. 
	Let $\varphi_0 \in \mathcal B ( \Teich( \mathbb H ) ) \cap B(\mathbb H^\ast, E^\ast )$ and 
	$f_0$ be the developing map of $\varphi_0$. Moreover, we denote Ahlfors' quasiconformal reflection of 
	quasicircle $f_0 ( \hat {\mathbb R} )$ by $\Lambda$, and $K_{\Lambda}$ is the maximal dilatation of $\Lambda$. 
	We fix a positive number $\epsilon$ with $2 \epsilon K_{\Lambda} < 1$ .  Then, we consider Ahlfors--Weil section
		\[
			\AW_{\varphi_0 , \Omega^\ast}: \Delta_{B(\Omega^\ast)} ( 0 ; \epsilon ) \to  \Bel ( \Omega ), 
		\]
	which is the section of 
		\[
			\beta_{\varphi_0} : \Bel ( \mathbb H ) \ni \mu \mapsto S ( f_\mu \circ f_0^{-1} ) \in B(\Omega^\ast)
		\]
	around the origin, where $\Omega : = f_0 ( \mathbb H ), \Omega^\ast : = f_{0}( \mathbb H^\ast )$. 
	Set $g: = f_\mu \circ f_0^{-1}$ on $\Omega^\ast$, 
	the Beltrami differntial $ \mu_g$ of a quasiconformal extenation of $g$ is  
		\[
			\widetilde{ \AW }_{\varphi_0 , \Omega^\ast} ( S(g) ) (w) = 
			\frac {\frac{1}{2} ( w - \Lambda ( w ) )^2 
			\cdot S(g) ( \Lambda ( w ) ) 
			\cdot \Lambda_{\bar w} ( w )}
			{1 + \frac{1}{2} ( w - \Lambda ( w ) )^2 
			\cdot S(g) ( \Lambda ( w ) ) 
			\cdot \Lambda_{w} ( w )} \ \ \ ( w \in \Omega), 
		\]
	where $w = f_0 ( \zeta )$. 
	From the definition of $\epsilon$, the $L^\infty$ norm of $\mu_g$ is less than $1$. 
	We will prove $\mu_g \in  \Bel (\Omega, f_0(E) )$. 
	To prove this, we observe the decay of $\mu_g$ when  $w$ moves away from the set  $f_0(E)$. From the chain rule 
	of the Shcwarizian derivative, we get
		\begin{align*}
			S(g) &= S ( f_\mu ) \circ f_0^{-1} \cdot ({f'}_0^{-1})^2 + S(f_0^{-1}),  \\
			S(f_0^{-1}) &= - S ( f_0 ) \circ f_0^{-1} \cdot ({f'}_0^{-1})^2 . 
		\end{align*}
	Moreover, since
		\[
			S(f_0^{-1}) ( w ) \rho^{-2}_{\Omega^\ast} ( w ) = - S ( f_0 ) ( \zeta ) \rho^{-2}_{\mathbb H ^\ast}(\zeta), 
		\]
	and $ f_0(E^\ast_{r/c}) \subset  f_0(E^\ast)_r$ ( Lemma \ref{lem: far from far} ),  for sufficiently large $r$, we get 
		\[
			\sup _{w \in \Omega^\ast \setminus ( f_0(E^\ast) )_r }
			| S(f_0^{-1}) ( w ) \cdot \rho^{-2}_{\Omega^\ast} ( w ) |
			\leq 
			\sup _{\mathbb H ^\ast \setminus ( E^\ast )_{r / c} } 
			| S ( f_0 ) ( \zeta ) \cdot \rho^{-2}_{\mathbb H ^\ast} ( \zeta) |. 
		\]
	Thus, $S(f_0^{-1}) \in B(\Omega^\ast,  f_0(E^\ast) )$. Similarly, we also get 
		\[
			 S ( f_\mu ) \circ f_0^{-1} \cdot ({f'}_0^{-1})^2 ( \zeta ) \cdot \rho^{-2}_{\Omega^\ast} ( \zeta ) 
			 = S ( f_\mu ) (z) \cdot \rho^{-2}_{\mathbb H^\ast} ( z). 
		\]
	Therefore, 
	if  $S(f_\mu) \in B(\mathbb H^\ast, E^\ast )$, 
	$S(g) \in B (\Omega^\ast, f_0(E^\ast) )$. 
	From properties of Ahlfors' quasiconformal reflection (Theorem \ref{Ahlfors ref}), $\mu_g$ satisfies
		\[
			| \mu_g ( w ) | \leq 
			\frac {K_{\Lambda} \cdot S(g)( \Lambda ( w ) ) \cdot \rho^-2_{\Omega^\ast} ( \Lambda ( w ) ) }
			{1 -  K_{\Lambda} \cdot S(g)(  \Lambda ( w ) ) \cdot \rho^-2_{\Omega^\ast} (  \Lambda ( w ) ) } 
		\]
	so, $\mu_g \in \Bel (\Omega, f_0(E) ) $. That is, we have proven  
		\[
			\widetilde{ \AW }_{\varphi_0, \Omega^\ast} : \Delta_{B(\Omega^\ast)} ( 0 ; \epsilon ) 
			\cap B(\Omega^\ast, f_0(E^\ast) )
			\to \Bel (\Omega, f_0(E) ). 
		\]
	Finally, we take a pull-back of them using Lemmas \ref{lem: pull of quad. diff} and Lemma \ref{lem: pull of bel. diff}.
	\qed
\end{proofbar}
\end{pr}

\section{Continued proof of Theorem \ref{main theorem}}

\begin{claim} \label{last claim}
\begin{lembar}	
	
	Under the setting of the proof of Theorem \ref{main theorem},one of the following is satisfied.
		\begin{align*}
			 \sup_{z \in \mathbb H^\ast} | ( \tilde \varphi ( z ) -\tilde \psi ( z ) )
			 \rho^{-2}_{\mathbb H^\ast} ( z ) | & \geq 1/2 \\
			\limsup_{r \to \infty}
			\sup_{z \in \mathbb H^\ast \setminus ( \tilde R_1^\ast )_{r} } 
			| \widetilde \psi ( z ) \rho^{-2}_{\mathbb H^\ast}(z) |
			&\geq 1/16
		\end{align*}
\end{lembar}	
\end{claim}

\begin{pr}
\begin{proofbar}

	We will prove this using the method of contradiction. That is, we will assume the negation of both: 
		\begin{align}
			 \sup_{z \in \mathbb H^\ast} | ( \tilde \varphi ( z ) -\tilde \psi ( z ) )
			 \rho^{-2}_{\mathbb H^\ast} ( z ) | & \leq 1/2 \label{first assume},  \\
			\limsup_{r \to \infty}
			\sup_{z \in \mathbb H^\ast \setminus ( \tilde R_1^\ast )_{r} } 
			| \widetilde \psi ( z ) \rho^{-2}_{\mathbb H^\ast}(z) |
			&\leq 1/16 \label{second assume}. 
		\end{align}
	We take $\mu \in \Bel ( \mathbb H, E )$ which is the Beltrami coefficient of the quasiconformal extension of 
	the developing map of $\tilde \varphi$. 
	Using $\mu$,  
	from \eqref{first assume} and Theorem\ref{Thm: translation map betweem Teichsp}, 
	it follows that there exists $\Psi \in B(\mathbb H^\ast, F^\mu )$
	such that $\tilde \psi= \mathcal B \circ a_{\mu} \circ ( \mathcal B^{\mu})^{-1} (\Psi) $. 
	We define $\hat { \Psi }: = ( j \Psi j ) \circ w^\mu \cdot ( w'^\mu )^2$, then
	$\hat {\nu}: =  \rho_{\Omega}^{-2} \overline {\hat { \Psi} }$ is a Beltrami coefficient of the quasiconformal extension of
	the developing map of $S(f_{\tilde \psi} \circ f_{\tilde \varphi} ^{-1} )$. 
	Moreover, the quasiconformal extension of the developing map of $\tilde \psi$ is given by
		
		\[
			\tau : = \frac{\hat {\nu} + \theta_\mu \overline{  \hat {\nu} \circ f_{\tilde \psi} } \cdot \mu}
			{1 + \theta_\mu \overline { \overline{  \hat {\nu} \circ f_{\tilde \psi} } \cdot \mu}  }, \ \ \
			 \text{where} \ \ \theta_\mu : = \frac {( f_{\tilde \varphi })_z } {\overline{ ( f_{\tilde \varphi })_z } }. 
		\]
	
	Since the absolute value  $| \hat { \nu }$ | is less than $1 / 4$ on $\Omega \setminus ( f_{\tilde \varphi} ( E ) )_r$ ,
	for sufficiently large $r$, from \eqref{second assume} and Remark\ref{rem: gardiner}, 
	$| \tau |$ satisfies $| \tau | \leq 3/4$ on $\mathbb H \setminus E_r$. 
	That is, $f_{\tilde \psi}$ is at most 7-quasiconformal map. 
	Thus by taking sufficiently large $r$ if necessary, from  Wolpert-Fujikawa's inequality (Theorem\ref{wol. fuj. ine}), 
	the length of a simple closed geodesic only be deformed by a factor of at most $7$.
	On the other hand, the quasiconformal map corresponding to $f_{\tilde \psi}$ is deformed by more than
	$100$ times, according to \eqref{lengths on tilde psi surface} .
	
	Here, if $( L_j )$ is contained in the region where $f_{\tilde \psi}$ can only be deformed by at most $7$ times,
	then it is immediately contradictory. Otherwise, consider the case where this is not the case.
	In this case, since $( L_j )$ gives a pants decomposition of the surface, there exists a simple closed
	geodesic $L$ that is contained in the region where $f_{\tilde \psi}$  is contained in a region that can only be
	deformed by a factor of at most $7$, and there is a simple closed geodesic $L$ that intersects with
	some $L_j$. At this time, because $L_j$ is shortened, $L$ is lengthened (Theorem \ref{thecollorlem}). 
	Therefore, it is still a contradiction.
	
\end{proofbar}
\end{pr}

\begin{pr}[Continued proof of Theorem \ref{main theorem}]
\begin{proofbar}

	First, is \eqref {first estimates} is grater than $1 / 4$, It has already been proven. 
	Otherwise, from Claim \ref{last claim}, 
	 we get 
		\begin{align*}
			\limsup _{r \to \infty}
			\sup_{z \in \mathbb H^\ast \setminus ( \tilde R_1^\ast )_r )} | \widetilde \psi ( z ) \rho^{-2}_{\mathbb H^\ast}(z) | 
			& \geq 1/16 \\
			\sup_{z \in \mathbb H^\ast \setminus ( \tilde R_1^\ast )_r )} | \tilde \varphi ( z ) \rho^{-2}_{\mathbb H^\ast}(z) | 
			& \to 0 \label{Davidfbgroup}
		\end{align*}
	Therefore, taking $\varepsilon < 1/64$ , 
		\[
			\| \varphi - \psi \| _{B(R)} \geq 1/16. 
		\]
	\qed
\end{proofbar}
\end{pr}

\Addresses

\end{document}